\documentclass[12pt,english,a4paper]{amsart}
\usepackage[T1]{fontenc}
\usepackage[latin1]{inputenc}
\usepackage{geometry}
\usepackage{graphicx}

\makeatletter
 \theoremstyle{plain}    
 \newtheorem{thm}{Theorem}[section]
 \numberwithin{equation}{section} 
 \numberwithin{figure}{section} 
 \theoremstyle{plain}
 \theoremstyle{plain}    
 \newtheorem{prop}[thm]{Proposition} 
 \theoremstyle{definition}
  \newtheorem{example}[thm]{Example}
 \theoremstyle{plain}    
 \newtheorem{lem}[thm]{Lemma} 
 \theoremstyle{remark}
 \newtheorem{rem}[thm]{Remark}
 \theoremstyle{plain}    
 \newtheorem{cor}[thm]{Corollary} 

\def\makebbb#1{
    \expandafter\gdef\csname#1\endcsname{
        \ensuremath{\Bbb{#1}}}
}
\makebbb{R}
\makebbb{N}
\makebbb{Z}
\makebbb{C}
\makebbb{G}
\makebbb{E}
\makebbb{H}
\makebbb{P}
\makebbb{B}
\makebbb{K}
\makebbb{Q}

\usepackage{babel}
\makeatother
\begin{document}

\title{Bergman kernels for weighted polynomials and weighted equilibrium
measures of $\C^{n}$}

\author{Robert Berman}

\email{robertb@math.chalmers.se}

\begin{abstract}
Various convergence results for the Bergman kernel of the Hilbert
space of all polynomials in $\C^{n}$ of total degree at most $k,$
equipped with a weighted norm, are obtained. The weight function $\phi$
is assumed to be $\mathcal{C}^{1,1},$ i.e. $\phi$ is differentiable
and all of its first partial derivatives are locally Lipshitz continuous.
The convergence is studied in the large $k$ limit and is expressed
in terms of the global equilibrium potential associated to the weight
function, as well as in terms of the Monge-Ampere measure of the weight
function itself on a certain set. A setting of polynomials associated
to a given Newton polytope, scaled by $k,$ is also considered. These
results apply directly to the study of the distribution of zeroes
of random polynomials and of the eigenvalues of random normal matrices.

\tableofcontents{}
\end{abstract}
\maketitle

\section{Introduction}

Let $\phi$ be a given function on $\C^{n}$ that will be refered
to as the \emph{weight function.} We will assume that $\phi$ is in
the class $\mathcal{C}^{1,1},$ i.e. $\phi$ is differentiable and
all of its first partial derivatives are locally Lipshitz continuous.
In particular, the {}``curvature form'' $dd^{c}\phi(z)$ (where
$d^{c}:=i(-\partial+\overline{\partial})/4\pi,$ so that $dd^{c}=i\partial\overline{\partial}/2\pi)$
is well-defined for almost any fixed point $z$ and its determinant
$\det(dd^{c}\phi):=(dd^{c}\phi)^{n}/\omega_{n}$ is in $L_{loc}^{1}(\C^{n}),$
where $\omega_{n}=(2\pi dd^{c}\left|z\right|^{2})^{n}/n!,$ $\omega_{n}$
is the Lebesgue measure on $\C^{n}.$ It is further assumed that $\phi$
is sufficiently large at infinity so as to make the following weighted
norms finite (assumption \ref{eq:ass on growth of phi} below). Denote
by $\mathcal{H}_{k}$ the Hilbert space of all polynomials $f_{k}$
of total degree less than $k$ in $\C^{n}$ equipped with the weighted
norm \[
\left\Vert f_{k}\right\Vert _{k\phi}^{2}:=\int_{\C^{n}}\left|f_{k}(z)\right|^{2}e^{-k\phi}\omega_{n},\]
 The \emph{Bergman kernel} of the Hilbert space $\mathcal{H}_{k}$
will be denoted by $K_{k}(z,w)$ (formula \ref{eq: def of K}). In
the final section of the paper we will also consider a variant of
this setting where $\mathcal{H}_{k}$ is replaced by the space of
all weighted polynomials associated to a given Newton polytope, scaled
by $k.$

Given this setup, we will consider three natural positive measures
on $\C^{n}.$ First the weighted \emph{equilibrium measure}\[
(dd^{c}\phi_{e})^{n}/n!,\]
 where \emph{}$\phi_{e}$ is the upper envelope \ref{eq:extem metric}.
Then the weak large $k$ limit of the following two sequence of measures
defined in terms of the Bergman kernel of the Hilbert space $\mathcal{H}_{k}:$
\begin{equation}
k^{-n}B_{k}(z)\omega_{n},\label{eq:intro B meas}\end{equation}
 where $B_{k}(z):=K_{k}(z,z)e^{-k\phi}$ will be referred to as the
\emph{Bergman function}%
\footnote{$B_{k}(z)$ corresponds, in the classical theory of orthogonal polynomials,
to the so called Christoffel function (see \cite{b-l} and references
therein). %
} \emph{}and of\[
(dd^{c}(k^{-1}\textrm{ln\,$K^{k}(z,z)))^{n}/n!$.}\]
It is not hard to see that the total integrals of all three measures
coincide. The main point of the present paper is to show the corresponding
\emph{local} statement. In fact, all three measures will be shown
to coincide with the measure \[
1_{D}(dd^{c}\phi)^{n}/n!\]
where $1_{D}$ is the characteristic function of a certain compact
subset $D$ of $\C^{n}$ (corollary \ref{cor:equil meas}, theorem
\ref{thm:B in L1} and \ref{thm:ln K}). The basic case is when $\phi(z)=\left|z\right|^{2}$.
Then the limiting measure is $\pi^{-n}$ times the Lebesgue measure
times the characteristic function of the unit ball. It should also
be pointed out that even stronger convergence results will be obtained,
where the regularity assumption on $\phi$ is crucial. 

Moreover, for any interior point of $D$ where $\phi$ is smooth and
where $dd^{c}\phi>0$ we will show (theorem \ref{thm:asymp expansion})
that the Bergman kernel $K_{k}(z,w),$ i.e. with two arguments, admits
a complete local asymptotic expansion in powers of $k,$ such that
the coefficients of the corresponding symbol expansion coincide with
the Tian-Zelditch-Catlin expansion for a positive Hermitian holomorphic
line bundle. Furthermore, it will be shown (theorem \ref{thm:k as meas})
that globally on $\C^{n}$\[
\begin{array}{lr}
k^{-n}\left|K_{k}(z,w)\right|^{2}e^{-k\phi(z)}e^{-k\phi(w)}\omega_{n}(z)\wedge\omega_{n}(w)\rightarrow1_{D}\Delta\wedge(dd^{c}\phi)^{n}/n!\end{array},\]
weakly as measures on $\C^{n}\times\C^{n}$, where $\Delta$ is the
current of integration along the diagonal in $\C^{n}\times\C^{n}.$

\subsection{Relations to random eigenvalues and zeroes}

In a companion paper the present results will be applied to the study
of various random processes. For example, consider the case of the
complex plane (i.e. $n=1).$ On one hand it is well-known that the
measure $B_{k}(z)\omega_{n}$ represents the expected distribution
of \emph{eigenvalues} of a random normal matrix, when the size $N(=k+1)$
tends to infinity (given a certain probability measure, induced by
$\phi,$ on the space of all $N\times N$ normal matrices) \cite{e-f,e-m}.
On the other hand, the measure $(dd^{c}(\textrm{ln\,$K^{k}(z,z)))^{n}/n!$}$
represents the expected distribution of \emph{zeroes} of a random
polynomial considered as an element of the Hilbert space $\mathcal{H}_{k},$
equipped with the induced Gaussian probability measure \cite{sz2,sz3}.
From this point of view the present results say that, in the corresponding
limits, the eigenvalues and the zeroes tend to behave as interacting
charges in the presence of the external potential field $\phi$ (compare
\cite{le}). Furthermore, in these limits they all become uniformly
distributed (with respect to the measure $(dd^{c}\phi)^{n}/n!)$ on
the set $D.$

\subsection{Comparison with previous results}

It should be emphasized that in the case when $n=1$ the weak convergence
for $k^{-n}B_{k}\omega_{n}$ in theorem \ref{thm:B in L1} and the
subsequent corollary was first obtained recently by Hedenmalm-Makarov
\cite{h-m}, generalizing a previous result of Elbau-Felder \cite{e-f}
(see also the paper \cite{e-m} of Etingof-Ma). These papers were
stimulated by the study of normal random matrices and various diffusion-controlled
growth processes, by the physicists Zabrodin, Wiegmann et al (see
the survey \cite{z} and references therein). For example, in one
physical model (the Hele-Shaw cell) the set $D$ and its complement
correspond to two incompressible fluids on the plane of very different
viscosities, e.g. water and oil. In another model the set $D$ is
identified with a quantum Hall droplet of electrons. 

In the higher dimensional case there seem to be almost no previous
results concerning the weak convergence of $k^{-n}B_{k}\omega_{n}$
in the present context. However, the case of a smooth plurisubharmonic
weight $\phi$ is closely related to the Tian-Zelditch-Catlin expansion
of the Bergman kernel in the line bundle setting (see section \ref{sub:Comparison-with-the-line}
below). Moreover, the results concerning weighted polynomials associated
to a Newton polytope, considered in section \ref{sec:Weighted-polynomials-associated},
generalize results of Shiffman-Zelditch \cite{sz4} - compare section
\ref{sec:Weighted-polynomials-associated} for the relation to their
results.

\subsubsection{Comparison with the ''large deviation'' approach }

The approach in \cite{h-m,e-f,e-m} follows the paper \cite{j} by
Johansson, where the \emph{real} case was treated (corresponding to
Hermitian matrices). Even though the present approach is completely
different, it may be illuminating to recall the approach in the cited
papers as it makes the physical interpretations more transparent.
The starting point is that, by a classical formula due to Heine in
the theory of orthogonal polynomials, the measure $B_{k}(z)\omega_{n}$
in \ref{eq:intro B meas} paired with a given test function $g,$
may be expressed as \begin{equation}
\int_{\C}B_{k}(z)g(z)\omega_{n}=\sum_{i}\int_{\C^{N}}g(z_{i})e^{-N^{2}E_{N}(z_{1},...,z_{N})}\omega_{n}(z_{1})\cdots\omega_{n}(z_{N})/Z^{N}.\label{eq:heine}\end{equation}
 Here $N=k+1$ and $E_{N}(z_{1},...,z_{N})$ is the {}``energy''
\[
E_{N}(z_{1},...,z_{N})=-\sum\sum_{i\neq j}\ln\left|z_{i}-z_{j}\right|N^{-2}+\sum_{i}\phi(z_{i})N^{-1}\]
of a system of $N$ charges, each with charge $N^{-1}$ situated at
the positions $z_{i}.$ The normalizing number $Z_{N}$ (corresponding
to $g=1)$ is the {}``partition function'' of the system. Following
a large deviation estimate in \cite{j}, it is shown that the leading
contribution to the integral over $\C^{N}$ in formula \ref{eq:heine},
comes from a sequence of minimizers $(z_{1},...,z_{N})_{*}$ of the
energy $E_{N}.$ Moreover, \[
N^{-1}\sum_{i}\delta_{z_{i*}}\rightarrow\mu,\]
 where $\mu$ is the unique minimizer of the functional \[
E_{\phi}(\mu):=-\int_{\C}\int_{\C}\ln\left|z-w\right|d\mu(z)\wedge d\mu(w)+\int_{\C}\phi(z)d\mu(z)\]
on the space of all positive Borel measures with total integral equal
to one. This gives\[
k^{-1}\int_{\C}B_{k}(z)g(z)\omega_{n}\rightarrow\int_{\C}gd\mu\]
Finally, the variational inequalities corresponding to the latter
functional $E_{\phi}$ give that $\mu$ may be obtained as $dd^{c}\phi_{e}$
where $\phi_{e}$ is the upper envelope \ref{eq:extem metric} (compare
the book \cite{s-t}). 

However, when $n>1$ there seems to be no (useful) explicit formula
playing the role of Heine's formula \ref{eq:heine}, which is related
to the fact that the Monge-Ampere operator which plays the role of
the Laplace operator is \emph{non-linear.} Instead the approach in
the present paper is to obtain the asymptotics of $B_{k}$ by a completely
different method based on local holomorphic Morse inequalities (compare
\cite{berm1,berm2}), combined with some pluripotential theory (see
\cite{ko,kl,de3}). For example, theorem \ref{thm:ln K} can be seen
as an $L^{2}-$version of a seminal result of Siciak in pluripotential
theory \cite{kl}. It was Siciak who first defined and studied the
function $\phi_{e},$ given a possibly non-smooth function $\phi,$
in the context of polynomial approximation and interpolation theory
\cite{kl}. For a very recent study of weighted pluripotential theory
and its relation to the more studied unweighted theory see \cite{b}.

\subsection{\label{sub:Comparison-with-the-line}Comparison with the line bundle
setting}

Some time after the first version of the preprint of the present paper
had appeared the author managed to obtain variants of the present
results in the setting of an arbitrary given line bundle $L\rightarrow X$
over a compact projective complex manifold $X$ \cite{berm4}. In
this latter setting the role of $\phi$ is played by a metric on $L.$
In the special case when $\phi$ has positive curvature the Bergman
kernel asymptotics then follow from the Tian-Zelditch-Catlin expansion
(compare \cite{sz2} and reference therein). The present setting is
closely related to the case when $X$ is complex projective space
$\P^{n}.$ However, the difference is that, due to the growth assumption
\ref{eq:ass on growth of phi}, $\phi$ does not extend as a \emph{smooth}
(or even locally bounded) metric over the hyperplane {}``at infinity'',
$\P^{n}-\C^{n},$ which is assumed in \cite{berm4}. In \cite{berm4}
it was also observed that the convergence in theorem \ref{thm:B in L1}
can be made slightly stronger using $\bar{\partial}$-estimates (compare
remark \ref{rem:ae}). But the author has decided to keep the original
approach in this paper due to its simplicity and since it may perhaps
be useful in other contexts. Finally, in the present paper it is verified
that the assumption on the $\mathcal{C}^{2}-$smoothnes on $\phi$
appearing in \cite{berm4} may be relaxed to assuming that $\phi$
is in the class $\mathcal{C}^{1,1}.$ The motivation for this improvement
is first of all that this is the natural class in view of the regularity
proposition \ref{pro:ref of equil pot} - the point is that even if
$\phi$ is assumed to be $\mathcal{C}^{\infty}-$smooth, the equilibrium
weight $\phi_{e}$ will usually only be $\mathcal{C}^{1,1}.$ Moreover,
the weaker regularity assumption might turn out to be useful when
trying to, by some approximation scheme, use the present results to
prove the following more general conjectures.

\subsection{Conjectures concerning more general cases}

It seems natural to consider the weak convergence of $k^{-n}B_{k}(z)\omega_{n}$
and $(dd^{c}(k^{-1}\textrm{ln\,$K^{k}(z,z)))^{n}/n!$}$ for a weight
function $\phi$ which is only assumed to be lower semi-continuous.
It may even be allowed to be equal to $-\infty$ outside a subset
$E,$ but then $\omega_{n}$ has to be replaced with a suitable measure
$\nu$ supported on $E.$

For example, in the case that $E=\C^{n}$ it seems reasonable to conjecture
that the corresponding weak convergence still holds. The issue is
the convergence of the Bergman measures \ref{eq:intro B meas}. Indeed,
the weak convergence of $(dd^{c}(k^{-1}\textrm{ln\,$K^{k}(z,z)))^{n}/n!$}$
does hold as is shown by a simple modification of the proof given
here (for the modification in the case of a general set $E$ see \cite{berm3}).
In the case where $E$ is compact this was shown by Bloom-Shiffman
\cite{bl-sh} in the unweighted case $(\phi=0)$ and Bloom \cite{b2}
in the weighted case using a different aproach. For even more general
conjectures, where one assumes that the measure $\nu$ satisfies a
Bernstein-Markov inequality, see \cite{berm3,b-l3} - the planar case
$\C$ was obtained very recently in \cite{b-l}. 

However, one point of the present paper is also to show that under
the assumption that $\phi$ be $\mathcal{C}^{1,1}$ stronger convergence
results hold: the Monge-Ampere measure $(dd^{c}\phi_{e})^{n}/n!$
has a density $\det(dd^{c}\phi_{e})$ that is in $L^{1}(\C^{n})$
and $k^{-n}B_{k}$ converges to $\det(dd^{c}\phi_{e})$ in $L^{1}(\C^{n})$.
In fact, one even has point-wise convergence almost everywhere on
$\C^{n}$ (see remark \ref{rem:ae}). Moreover, the rate of convergence
of $k^{-1}\textrm{ln\,$K^{k}(z,z)$}$ towards $\phi_{e}(z)$ is optimal.

\section{Weighted equilibrium measures}

Given a weight function $\phi$ on $\C^{n}$ the corresponding {}``equilibrium
potential'' or {}``equilibrium weight'' $\phi_{e}$ is defined
by

\begin{equation}
\phi_{e}(z)=\sup\left\{ \widetilde{\phi}(z):\,\,\,\widetilde{\phi}\in\mathcal{L}:\,\widetilde{\phi}\leq\phi\,\,\textrm{on\,$\C^{n}$}\right\} .\label{eq:extem metric}\end{equation}
 where $\mathcal{L}(\C^{n})$ is the \emph{Lelong class}, consisting
of all plurisubharmonic functions in $\C^{n}$ of (at most) logarithmic
growth at infinity: \begin{equation}
\mathcal{L}(\C^{n}):=\left\{ \psi:\psi\,\textrm{psh,\,\,}\psi(z)\leq\ln^{+}\left|z\right|^{2}+O(1)\right\} \label{eq:def of lelong}\end{equation}
(note that we are using a non-standard exponent $2$ in the growth
assumption for the Lelong class).

Then $\phi_{e}$ is also in $\mathcal{L}(\C^{n})$ (compare the appendix
in \cite{s-t}). The corresponding \emph{weighted equilibrium measure}
is the Monge-Ampere measure $(dd^{c}\phi_{e})^{n}/n!$ (see \cite{ko}
for the definition of the Monge-Ampere measure of a locally bounded,
possibly non-smooth, plurisubharmonic function). 

Consider the set\begin{equation}
D:=\{\phi_{e}=\phi\}\subset\C^{n},\label{eq:def of D}\end{equation}
which is closed since $\phi_{e}$ is usc. We will assume that \begin{equation}
\phi(z)\geq(1+\epsilon)\ln\left|z\right|^{2},\,\,\textrm{when\,}\left|z\right|>>1\label{eq:ass on growth of phi}\end{equation}
 for some positive number $\epsilon,$ to make sure the $L^{2}-$norms
in section \ref{sec:Bergman-kernels-for} are finite. Then the set
$D$ is bounded (since, by definition, $\phi_{e}$ has logarithmic
growth, corresponding to $\epsilon=0$ above) and hence compact. 

The following regularity result plays a crucial role in the proof
of theorem \ref{thm:B in L1} below. Its proof follows closely the
exposition in \cite{de3} of the approach of Bedford-Taylor for the
Dirichlet problem in the unit-ball for the Monge-Ampere equation \cite{b-t}.
See also \cite{b} (in particular theorem 3) for a similar situation,
where the relation to \emph{free} boundary value problems for the
Monge-Ampere equation is pointed out. The main new feature here is
that growth condition at {}``infinity'' has to be taken into account.

\begin{prop}
\label{pro:ref of equil pot}Suppose that $\phi$ is a function on
$\C^{n}$ in the class $\mathcal{C}^{1,1},$ satisfying the growth
assumption \ref{eq:ass on growth of phi}. Then $\phi_{e}$ is also
in the class $\mathcal{C}^{1,1}.$ In particular, the Monge-Ampere
measure of $\phi_{e}$ is absolutely continuous w.r.t Lesbegue measure
and coincides with the corresponding $L_{loc}^{\infty}$ $(n,n)-$form
obtained by a point-wise calculation: \begin{equation}
(dd^{c}\phi_{e})^{n}=\det(dd^{c}\phi_{e})\omega_{n}\label{eq:ptwise repr of equil meas}\end{equation}
Moreover, the following identity holds almost everywhere on the set
$D=\{\phi_{e}=\phi\}:$ \begin{equation}
\det(dd^{c}\phi_{e})=\det(dd^{c}\phi)\label{eq:monge on D}\end{equation}

\end{prop}
\begin{proof}
To obtain the $\mathcal{C}^{1,1}-$regularity first observe that there
is a positive number $R$ such that \begin{equation}
\left|z\right|>R\Rightarrow\,\,\,\phi_{e}(z+h)\leq\phi(z)+2R^{-1}h\label{eq:log est}\end{equation}
Indeed, since $\phi_{e}(z)$ is in the Lelong class there is a constant
$C$ such that $\phi_{e}(z+h)\leq\textrm{$2($ln $\left|z+h\right|+C\leq\textrm{$2($ln $(\left|z\right|)+C)+2\left|z\right|^{-1}\left|h\right|$)$\leq\phi(z)+2R^{-1}h$ }$}$for
$\left|z\right|$ sufficintly large, using the growth assumption on
$\phi$ in the last step. 

\emph{Step1:} \textbf{\emph{}}\emph{$\phi_{e}$ is Lipschitz continuous}

First assume that $\left|z\right|<R.$ Then \begin{equation}
\phi_{e}(z+h)\leq\phi(z+h)\leq\phi(z)+(\sup_{\left|z\right|\leq R}\left|d\phi\right|)h,\label{eq:lip est}\end{equation}
 using that $\phi$ is assumed to be $\mathcal{C}^{1}-$smooth in
the last step. Combining \ref{eq:log est} and \ref{eq:lip est} now
gives that $\phi_{e}(z+h)-C_{R}h$ with $C_{R}=\max\{2R^{-1},(\sup_{\left|z\right|\leq R}\left|d\phi\right|)\}$
is a candidate for the sup defining $\phi_{e}$ and hence bounded
from above by $\phi_{e}(z),$ i.e. \[
\phi_{e}(z+h)-\phi_{e}(z)\leq C_{R}h.\]
 Applying the previous inequality with $z$ and $h$ replaced by $z+h$
and $-h,$ respectively, finally finishes the proof of step 1.

\emph{Step2:} \textbf{\emph{}}\emph{$d\phi_{e}$ exists and is Lipschitz
continuous}

Following the exposition in \cite{de3} of the approach of Bedford-Taylor
it is enough to prove the following inequality:

\begin{equation}
(\phi_{e}(z+h)+\phi_{e}(z-h))/2-\phi(z)\leq C\left|h\right|^{2},\label{eq:pf of thm reg ineq 2}\end{equation}
where the constant only depends on the second derivatives of $\phi.$
Indeed, given this inequality (combined with the fact that $\phi_{e}$
is psh) a Taylor expansion of degree $2$ gives the following bound
close to a fixed point $z_{0}$ for a local smooth approximation $\phi_{\epsilon}$
of $\phi_{e}:$ \[
\left|D^{2}\phi_{\epsilon}\right|\leq C\]
 where $\phi_{\epsilon}:=\phi_{e}*u_{\epsilon},$ using a a local
regularizing kernel $u_{\epsilon}$ and where $D^{2}\phi_{\epsilon}$
denotes the real local Hessian matrix of $\phi_{\epsilon}.$ Letting
$\epsilon$ tend to $0$ then proves Step 2. Finally, to see that
the inequality \ref{eq:pf of thm reg ineq 2} holds we repeat the
argument in Step 1 after replacing $\phi_{e}(z+h)$ by the following
element of the Lelong class: \[
g(z):=(\phi_{e}(z+h)+\phi_{e}(z-h))/2\]
 to get $g(z):=(\phi_{e}(z+h)+\phi_{e}(z-h))/2\leq(\phi(z+h)+\phi(z-h))/2\leq C_{R}\left|h\right|^{2}$
now using that $\phi$ is assumed to be $\mathcal{C}^{1,1}$ in the
last step (by combining the estimate \ref{eq:lip est} with a simple
variant of the estimate \ref{eq:log est} for $\left|z\right|>R).$
Hence, as above we conclude that $g(z)-C_{R}h^{2}\leq\phi_{e}(z)$
which proves \ref{eq:pf of thm reg ineq 2} and hence finishes the
proof of step 2. 

Now, by the $\mathcal{C}^{1,1}-$regularity, the derivatives $\frac{\partial^{2}\phi}{\partial z_{i}\partial\bar{z_{j}}}\phi_{e}$
are in $L_{loc}^{\infty}$ and it is well-known that this implies
the identity \ref{eq:ptwise repr of equil meas} for the Monge-Ampere
measure \cite{de3}. Finally, to see that \ref{eq:monge on D} holds,
it is enough to prove that \[
\frac{\partial^{2}\phi}{\partial z_{i}\partial\bar{z_{j}}}(\phi_{e}-\phi)=0\]
almost everywhere on $D=\{\phi_{e}=\phi\}.$ To this end we apply
a calculus lemma in \cite{k-s} (page 53) to the $\mathcal{C}^{1,1}-$function
$\phi_{e}-\phi$ (following the approach in \cite{h-m}), which even
gives the corresponding identity between all real second order partial
derivatives.
\end{proof}
\begin{example}
\label{rem:C11 counter}Even if $\phi$ is assumed to be $\mathcal{C}^{\infty}-$smooth,
then $\phi_{e}$ will usually only be $\mathcal{C}^{1,1}.$ For example,
if $n=1$ and $\phi$ is radial and smooth then $\phi_{e}$ is generically
not $\mathcal{C}^{2}-$cf. example 5.1 in \cite{berm4}. Indeed, writing
$\phi(z)=\Phi(v)$ in terms of the logarithmic radial coordinate $v=\ln\left|z\right|^{2}$
the equilibrium potential $\phi_{e}$ corresponds to $\Phi_{e}$ obtained
as a convex envelope: $\Phi_{e}(v)=\sup\{\widetilde{\Phi}(v):\,\widetilde{\Phi}$
convex, $\left|\frac{d\widetilde{\Phi}}{dv}\right|\leq1,\,\widetilde{\Phi}\leq\Phi\,\textrm{on$\,\R$}\}.$
Hence, the graph of $\Phi_{e}$ may be obtained as a convex hull,
which will clearly {}``generically'' not be $\mathcal{C}^{2}$ (compare
figure \ref{cap:intro}). %
\begin{figure}
\begin{center}\includegraphics{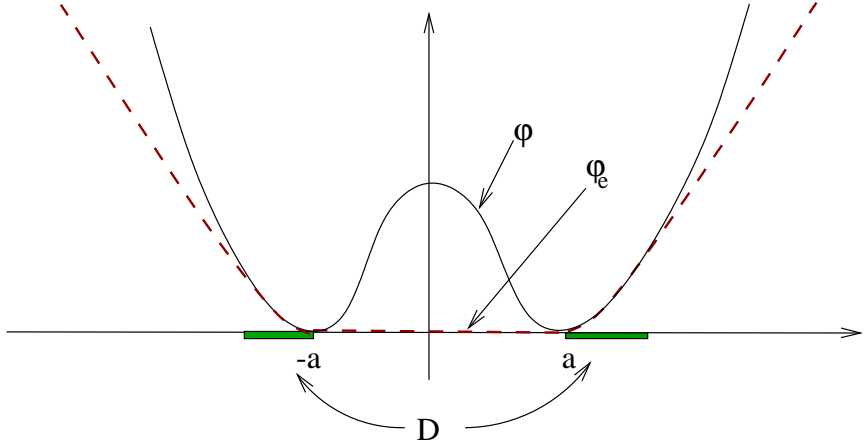}\end{center}

\caption{\label{cap:intro}A radial weight function $\phi$ on $\C_{z}$ may
be represented as a function on $\R_{v},$ $v=\ln\left|z\right|^{2}.$
Then the graph of the equilibrium potential $\phi_{e}$ (the dashed
line) is obtained as a \emph{convex hull} of the graph of $\phi$
(subject to an extra slope condition forcing $\left|\frac{d\phi_{e}}{dv}\right|\leq1).$
The set $D$ is then the projection of the set where the two graphs
coincide. In the case of the figure the outer boundary of $D$ is
given by the set where $\left|\frac{d\phi}{dv}\right|=1.$ In particular,
$\phi_{e}$ will note be twice differentiable at the two minima of
$\phi$ as long as they are non-degenerate, i.e. at the inner boundary
of the set $D,$ nor at the outer boundary of $D$ (generically).}
\end{figure}

\end{example}
Before turning to Bergman kernels we state the following {}``domination
lemma'' that we will have use for later. We refer to theorem 2.6
in the appendix by Bloom in \cite{s-t} where further references can
be found. 

\begin{lem}
\label{lem:domin}Suppose that $f$ is a polynomial of degree less
than $k$ in $\C^{n}$ such that \[
\left|f\right|^{2}e^{-k\phi}\leq C\,\,\,\textrm{on\,$D$.}\]
Then \[
\left|f\right|^{2}\leq Ce^{k\phi_{e}}\,\,\,\textrm{on\,$\C^{n}$.}\]

\end{lem}
The previous lemma is a simple consequence of the following {}``domination
principle'' applied to $\phi_{1}=\ln(\left|f\right|^{2}/C),\,\phi_{2}=\phi_{e}.$
Suppose that $\phi_{1},\phi_{2}\in\mathcal{L}(\C^{n})$ and that $\phi_{2}$
is locally bounded and of maximal growth (i.e. $\phi_{2}(z)=\ln^{+}\left|z\right|^{2}+O(1)$).
Then $\phi_{1}\leq\phi_{2}$ on the support of $(dd^{c}\phi_{e})^{n}$
implies $\phi_{1}\leq\phi_{2}$ everywhere. This principle in turn
may be deduced from the following monotonicity property of the global
Monge-Ampere mass (lemma 15.4 in \cite{de3}):\begin{equation}
\phi_{1}\leq\phi_{2}+C\,\Rightarrow\int_{\C^{n}}(dd^{c}\phi_{1})^{n}\leq\int_{\C^{n}}(dd^{c}\phi_{2})^{n},\label{eq:comp of mass}\end{equation}
for any $\phi_{1},\phi_{2}\in\mathcal{L}(\C^{n})$ (see the appendix
in \cite{b-b} for a considerably more general context). The proof
of \ref{eq:comp of mass} in \cite{de3} just uses that $\mathcal{L}(\C^{n})$
is closed under taking the $\max$ of two elements. One final consequence
of \ref{eq:comp of mass} that we shall have use for is that \begin{equation}
\int_{\C^{n}}(dd^{c}\phi_{1})^{n}\leq1\label{eq:mass bound in lelong}\end{equation}
 for any $\phi_{1}$ in $\mathcal{L}(\C^{n})$ (just take $\phi_{2}(z)=\ln^{+}\left|z\right|^{2}).$

\begin{rem}
\label{rem:The-growth-assumption}The growth assumption \ref{eq:ass on growth of phi}
is used to make sure that the weighted $L^{2}-$norms in the following
section are finite. However, the purely pluripotential results in
this section hold, with the same proofs, under the slightly weaker
assumption that $\phi(z)-\ln\left|z\right|^{2}\rightarrow\infty$
as $\left|z\right|\rightarrow\infty$ (which is the growth assumption
appearing in \cite{s-t}). In particular, $D$ is still compact then.
There is also a variant of the theory where one allows $\phi(z)-\ln\left|z\right|^{2}$
to be bounded (but then $D$ will in general be non-bouned). The point
is that then $\phi$ extends to a locally bounded metric on the hyperplane
line bundle $\mathcal{O}(1)\rightarrow\P^{n}$ over complex projective
space $\P^{n}.$ However, for the analog of proposition \ref{pro:ref of equil pot}
to hold one must further assume that $\phi$ extends to a \emph{smooth}
metric on $\mathcal{O}(1)\rightarrow\P^{n}.$ For the general line
bundle case see \cite{berm4}. 
\end{rem}

\section{\label{sec:Bergman-kernels-for}Bergman kernels for weighted polynomials}

Denote by $\mathcal{H}_{k}$ the Hilbert space of all polynomials
$f_{k}$ of total degree less than $k$ in $\C^{n}$ with the weighted
norm \[
\left\Vert f_{k}\right\Vert _{k\phi}^{2}:=\int_{\C^{n}}\left|f_{k}(z)\right|^{2}e^{-k\phi}\omega_{n}\]
(by the assumption \ref{eq:ass on growth of phi} the norm is finite
for $k$ sufficiently large). The \emph{Bergman kernel} of the Hilbert
space $\mathcal{H}_{k}$ may be defined as \begin{equation}
K_{k}(z,w):=\sum_{i}\psi_{i}(z)\overline{\psi_{i}(w)}.\label{eq: def of K}\end{equation}
 where $\psi_{i}$ is any given orthonormal bases of $\mathcal{H}_{k}.$%
\footnote{$K_{k}$ may also be defined as the reproducing kernel of $\mathcal{H}_{k}$
or as the integral kernel of the orthogonal projection from all smooth
functions on the subspace $\mathcal{H}_{k}.$%
} We will call \begin{equation}
B_{k}(z):=K_{k}(z,z)e^{-k\phi}=\sum_{i}\left|\psi_{i}(z)\right|^{2}e^{-k\phi}\label{eq:def of B}\end{equation}
the \emph{Bergman function} of $\mathcal{H}_{k}.$ It has the following,
very useful, extremal characterization: \begin{equation}
B_{k}(z)=\sup_{f_{k}\in\mathcal{H}_{k}}\left|f_{k}(z)\right|^{2}e^{-k\phi(z)}/\left\Vert f_{k}\right\Vert _{k\phi}^{2},\label{eq:extrm car of B}\end{equation}
 which follows from the reproducing property of the Bergman kernel.
Moreover, integrating \ref{eq:def of B} shows that $B_{k}$ is a
{}``dimensional density'' of the space $\mathcal{H}_{k}:$ \begin{equation}
\int_{\C^{n}}B_{k}\omega_{n}=\dim\mathcal{H}_{k}=1/n!k^{n}+o(k^{n})\label{eq:dim formel for B}\end{equation}
The following {}``local holomorphic Morse inequality'' estimates
$B_{k}$ point-wise from above. 

\begin{lem}
\emph{\label{lem:(Local-Morse-inequalities)}(Local Morse inequalities)
Suppose that the weight function $\phi$ is in the class $\mathcal{C}^{1,1}$.}
Given a compact $K$ set in $\C^{n}$ there is a constant $C_{K}$
such that \[
(i)\,\,\, k^{-n}B_{k}\leq C_{K}\]
 on $K.$ Moreover, the following bound holds at any point $z$ where
the second derivatives of $\phi$ exist (the complement of all such
points is a set which has zero measure w.r.t the Lesbegue measure):\[
(ii)\,\,\,\limsup_{k}k^{-n}B_{k}(z)\leq1_{X(0)}\det(dd^{c}\phi)(z),\]
 where $X(0)$ is the set where $dd^{c}\phi$ exists and is strictly
positive.
\end{lem}
See \cite{berm1} for the more general corresponding result for $\overline{\partial}-$harmonic
$(0,q)$- forms with values in a high power of an Hermitian line bundle.
Since we consider the slightly more general case (for $q=0)$ when
$\phi$ is merely in the class $\mathcal{C}^{1,1}$ we will go through
the proof and provide the necessery modifications.

\begin{proof}
First observe that we may without loss of generality subtract from
$\phi$ any function of the form $2\textrm{Re $g$ where }$$g$ is
an entire function. Indeed, the point-wise norm of a holomorphic function,
$\left|f(z)\right|^{2}e^{-k\phi(z)},$ is invariant under the transformation
$\phi\rightarrow\phi-2\textrm{Re $g$, }f\rightarrow f\textrm{exp$(kg)$}$
(geometrically this corresponds to a change of frame on the trivial
holomorphic line bundle). In particular, since $\phi$ is $\mathcal{C}^{1}$
we may (by taking $g(z)$ to be of the form $a+\sum_{i}a_{i}z_{i})$
assume that the first order jet of $\phi$ vanishes at any fixed point
which we after a change of coordinate may take to be zero: $\phi(0)=d\phi(0)=0.$
Similarly, if the second deriviatives of $\phi$ exist at $z=0$ we
may (after perhaps changing the coordinates by a unitary transformation)
and will assume that the second order jet of $\phi$ is given by $\sum_{i=1}^{n}\lambda_{i}\left|z_{i}\right|^{2}$
(by taking $g(z)$ to be of the form $\sum_{i,j}a_{i,j}z_{i}z_{j}).$
Now observe that the following holds in general: \begin{equation}
\left|\phi(z)\right|\leq C\left|z\right|^{2}\label{eq:lip est}\end{equation}
 Moreover, if the second derivatives of $\phi$ exist at $z=0,$ then
for any fixed positive number $R$ the following uniform convergence
holds when $k$ tends to infinity 

\begin{equation}
\sup_{\left|z'\right|\leq R}\left|k\phi(\frac{z'}{\sqrt{k}})-\sum_{i=1}^{n}\lambda_{i}\left|z_{i}'\right|^{2}\right|\rightarrow0\label{convergence of fiber m (new)}\end{equation}
To see this write $\psi(t):=\phi(tz)-\sum_{i=1}^{n}\lambda_{i}\left|tz_{i}\right|^{2}$
for $t\in[0,1].$ Since $\psi(0)=0$ and $\psi$ is $\mathcal{C}^{1}$
we have $\psi(1)=\int_{0}^{1}\psi'(t)dt,$ where $\psi'(t):=\frac{d\psi}{dt}(t).$
Moreover, by the Lipschitz assumption on the differential $d\phi$
(which vanishes at $z=0)$ we get $\left|\psi'(t)\right|\leq C\left|z\right|\left|z\right|$
which proves $(i).$ Morever, if the second derivatives of $\phi$
exist (and hence those of $\phi(z)-\sum_{i=1}^{n}\lambda_{i}\left|z_{i}\right|^{2}$
vanish) at $z=0,$ we get that for any given $\epsilon>0$ there is
a $\delta>0$ such that \[
\left|z\right|<\delta\Rightarrow\left|\psi'(t)\right|\leq\epsilon\left|z\right|\left|z\right|.\]
 In particular, if $\left|z'\right|\leq R$ we get $k\left|\phi(\frac{z'}{\sqrt{k}})-\sum_{i=1}^{n}\lambda_{i}\left|\frac{z_{i}'}{\sqrt{k}}\right|^{2}\right|\leq k\epsilon(\frac{R}{\sqrt{k}})^{2}=\epsilon R.$
But since $R$ is fixed and this holds for any $\epsilon>0$ the estimate
$(ii)$ follows. 

Now by the extremal characterization \ref{eq:extrm car of B} of the
Bergman function we get \[
k^{-n}B_{k}(0)\leq\frac{\left|f_{k}(0)\right|^{2}}{k^{n}\int_{\left|z\right|\leq R/\sqrt{k}}\left|f_{k}(z)\right|^{2}e^{-k\phi(z)}\omega_{n}}\]
By \ref{eq:lip est} we may replace $\phi(z)$ by $C\left|z\right|^{2}$
and then applying the submean property of holomorphic functions in
each polydisc $\{\left|z\right|=r\},$ $r\leq R/\sqrt{k},$ deduce
that \[
k^{-n}B_{k}(0)\leq\frac{\left|f_{k}(0)\right|^{2}}{k^{n}\int_{\left|z\right|\leq R/\sqrt{k}}\left|f_{k}(z)\right|^{2}e^{-Ck\left|z\right|^{2}}\omega_{n}}\leq\frac{1}{\int_{\left|z'\right|\leq R}e^{-C\left|z'\right|^{2}}\omega_{n}}<\infty\]
where we have also performed the change of variables $z'=z\sqrt{k}$
in the last integral. Finally, since $d\phi$ is locally Lipschitz
continuous it is not hard to check that, with say $R=1,$ the constant
$C$ may be taken to only depend on the compact set $K$ and not on
the fixed point. This proves $(i)$ in the statement of the lemma.
\end{proof}
Finally, to prove $(ii)$ note that a similar argument, gives \[
k^{-n}B_{k}(0)\leq\rho_{k,R}\frac{1}{\int_{\left|z'\right|\leq R}e^{-\sum_{i=1}^{n}\lambda_{i}\left|z_{i}'\right|^{2}}\omega_{n}},\]
where $\rho_{k,R}:=\textrm{exp}(\sup_{\left|z'\right|\leq R}\left|k\phi(\frac{z'}{\sqrt{k}})-\sum_{i=1}^{n}\lambda_{i}\left|z_{i}'\right|^{2}\right|),$
which by \ref{convergence of fiber m (new)} tends to one when $k$
tends to infinity, for any fixed positive number $R.$ Hence, \[
\limsup_{k}k^{-n}B_{k}(0)\leq\frac{1}{\int_{\left|z'\right|\leq R}e^{-\sum_{i=1}^{n}\lambda_{i}\left|z_{i}'\right|^{2}}\omega_{n}}\]
 and letting $R$ tend to infinity finally proves $(ii)$ in the statement
of the lemma, since the invers of the (normalized) Gaussian integral
equals $\lambda_{1}\lambda_{2}\cdot\cdot\cdot\lambda_{n}/(2\pi)^{n}$
if all eigenvalues $\lambda_{i}$ are positive and vanishes otherwise. 

In fact, the previous proof gives the following stronger variant of
$(ii):$ \begin{equation}
\limsup_{R}\limsup_{k}k^{-n}\left|f_{k}(z)\right|^{2}e^{-k\phi(z)}/\left\Vert f_{k}\right\Vert _{k\phi,\Delta_{k,R}}^{2}\leq1_{X(0)}(z)\det(dd^{c}\phi)(z),\label{eq:local morse on ball}\end{equation}
 where $f_{k}$ is holomorphic function defined in a fixed neighbourhood
of the point $z,$ where we have assumed that the second derivatives
of $\phi$ exist and where $\Delta_{k,R}$ denotes a polydisc centered
at $z$ of radius $R/\sqrt{k}.$ 

\begin{rem}
\label{rem:holder}The regularity assumption on $\phi$ is essentially
optimal for the previous lemma to hold. Indeed, if $\phi$ is only
assumed to be locally in the Hölder class $\mathcal{C}^{1,1-\delta}$
for $\delta\in]0,1]$ , then the uniform bound $(i)$ in the previous
lemma may not hold. To see this let $\phi(z)=\left|z\right|^{2-\delta}.$
Then \begin{equation}
B_{k}(0)\geq\left|1\right|^{2}/\int_{\C^{n}}\left|1\right|^{2}e^{-k\left|z\right|^{2-\delta}}\omega_{n}\geq Ck^{n/(1-\delta/2)},\label{eq:holder ex}\end{equation}
 where $C=\int_{\C^{n}}e^{-\left|z\right|^{2-\delta}}\omega_{n},$
contradicting $(i).$ 
\end{rem}
Using global information, the estimate in the previous lemma can be
considerably sharpened on the complement of $D$ (formula \ref{eq:def of D}),
as shown by the following lemma:

\begin{lem}
\label{lem:exponentia decay}The following inequality holds on all
of $\C^{n}$: \begin{equation}
B_{k}k^{-n}\leq C_{k}e^{-k(\phi-\phi_{e})}\label{eq:exp decay}\end{equation}
 where the sequence $C_{k}$ of positive numbers tends to $\sup_{z\in D\cap X(0)}\det(dd^{c}\phi)(z).$
In particular, there is a uniform constant $C$ such that \begin{equation}
B_{k}(z)k^{-n}\leq C\left|z\right|^{-2\epsilon k},\,\,\,\textrm{if\,\,$\left|z\right|>C.$}\label{eq:dom by L1}\end{equation}

\end{lem}
\begin{proof}
By the extremal property \ref{eq:extrm car of B} of $B_{k}$ it is
enough to prove the lemma with $B_{k}k^{-n}$ replaced by $\left|f_{k}\right|^{2}e^{-k\phi}$
for any element $f_{k}$ in $\mathcal{H}_{k}$ with global norm equal
to $k^{-n}.$ Since $D$ is compact the Morse inequalities in the
previous lemma give that\[
\left|f_{k}(z)\right|^{2}e^{-k\phi(z)}\leq C_{k},\,\,\, z\in D\]
 with $C_{k}$ as in the statement of the present lemma. Now the exponential
decay \ref{eq:exp decay} follows from lemma \ref{lem:domin}. Finally,
the uniform bound \ref{eq:dom by L1} is obtained from the growth
assumption \ref{eq:ass on growth of phi} on $\phi$.
\end{proof}
Now we are ready for the proof of the following main 

\begin{thm}
\label{thm:B in L1}Let $B_{k}$ be the Bergman function of the Hilbert
space $\mathcal{H}_{k}$ of weighted polynomials in $\C^{n}$ of degree
less than $k.$ Then \begin{equation}
k^{-n}B_{k}\rightarrow1_{D\cap X(0)}\det(dd^{c}\phi),\label{eq:l1 conv of B}\end{equation}
 in $L^{1}(\C^{n}).$ In particular, $k^{-n}B_{k}\omega_{n}$ conerges
weakly to the corresponding equilibrium measure $(dd^{c}\phi_{e})^{n}/n!.$
\end{thm}
\begin{proof}
First observe that, by the exponential decay in lemma \ref{lem:exponentia decay},
\[
\lim k^{-n}B_{k}(z)=0,\,\, z\in D^{c}\]
and, by \ref{eq:dom by L1} in the same lemma, the dominated convergence
theorem then gives \begin{equation}
\lim\int_{D^{c}}k^{-n}B_{k}\omega_{n}=0\label{pf of thm B: int B on D compl}\end{equation}
Next, observe that it is enough to prove that \begin{equation}
\lim\int_{D}k^{-n}B_{k}\omega_{n}=\int_{D\cap X(0)}(dd^{c}\phi)^{n}/n!\label{pf of thm B: claim}\end{equation}
Indeed, given this equality the local Morse inequalities (lemma \ref{lem:(Local-Morse-inequalities)}),
then force the convergence \ref{eq:l1 conv of B} on the compact set
$D.$ The proof proceeds precisely as in \cite{berm2} (part 1, section
2).

Finally, to prove that \ref{pf of thm B: claim} does hold, first
note that \[
1/n!=\lim\int_{\C^{n}}k^{-n}B_{k}\omega_{n}=\lim\int_{D}k^{-n}B_{k}\omega_{n}\]
using \ref{eq:dim formel for B} in the first equality and \ref{pf of thm B: int B on D compl}
in the second one. Now the local Morse inequalities (lemma \ref{lem:(Local-Morse-inequalities)})
applied to the compact set $D,$ give\[
1/n!=\lim\int_{D}k^{-n}B_{k}\omega_{n}\leq\int_{D\cap X(0)}(dd^{c}\phi)^{n}/n!\]
(where we have used the dominated convergence theorem on $D$). By
formula \ref{eq:monge on D} in proposition \ref{pro:ref of equil pot}
we may replace $\phi$ with $\phi_{e}$ in the right hand side, giving
\[
1/n!=\lim\int_{D}k^{-n}B_{k}\omega_{n}\leq\int_{D\cap X(0)}(dd^{c}\phi)^{n}/n!=\int_{D\cap X(0)}(dd^{c}\phi_{e})^{n}/n!\leq1/n!,\]
using \ref{eq:mass bound in lelong} in the last equality. But this
can only happen if all inequalities above are actually equalities,
which proves \ref{pf of thm B: claim} and finishes the proof of the
theorem, since it also shows that $1_{D\cap X(0)}(dd^{c}\phi)^{n}=(dd^{c}\phi_{e})^{n}$
as measures.
\end{proof}
\begin{rem}
\label{rem:ae}In fact, the slightly stronger statement that $k^{-n}B_{k}\rightarrow1_{D\cap X(0)}\det(dd^{c}\phi)$
point-wise almost everywhere on $\C^{n}$ holds. Moreover, the set
where the convergence fails may be described in terms of properties
of the {}``second order jets'' of $\phi$ and $\phi_{e}.$ This
is shown by using Hörmander's weighted $\overline{\partial}-$estimates
with weight $k\phi_{e}+\ln(1+\left|z\right|^{2}),$ precisely as in
\cite{berm4}. 
\end{rem}
\begin{cor}
\label{cor:equil meas}The weighted equilibrium measure corresponding
to the weight function $\phi$ is given by \[
(dd^{c}\phi_{e})^{n}/n!=1_{D\cap X(0)}(dd^{c}\phi)^{n}/n!=1_{D}(dd^{c}\phi)^{n}/n!\]
where $D=\{\phi_{e}=\phi\},$ which is compact under the assumption
\ref{eq:ass on growth of phi}. Moreover, the equality holds almost
everywhere w.r.t the Lesbegue measure $\omega_{n}$ (recall that by
proposition \ref{pro:ref of equil pot} $(dd^{c}\phi_{e})^{n}$ is
absolutely continuous w.r.t $\omega_{n})$. 
\end{cor}
\begin{proof}
The first equality in the corollary was shown in the end of the previous
proof. The final one follows from the fact that $dd^{c}\phi\geq0$
on $D$ (where it is defined); see Proposition 3.1 in \cite{berm4}.
\end{proof}
Next, we will show that the equilibrium potential $\phi_{e}$ may
be obtained from the logarithm of the Bergman kernel. This can be
seen as an $L^{2}-$version of a seminal result of Siciak on polynomial
approximation (see \cite{kl}).

\begin{thm}
L\label{thm:ln K}et $K_{k}$ be the Bergman kernel of the Hilbert
space $\mathcal{H}_{k}$ of weighted polynomials in $\C^{n}$ of degree
less than $k.$ Then
\end{thm}
\begin{equation}
k^{-1}\textrm{ln\,$K^{k}(z,z)\rightarrow\phi_{e}(z)$}\label{eq:conv of ln K in theorem ln K}\end{equation}
uniformly on $\C^{n}$ (the rate of convergence is of the order $n\ln k/k).$
In particular,

\begin{equation}
(dd^{c}(k^{-1}\textrm{ln\,$K^{k}(z,z)))^{n}\rightarrow(dd^{c}\phi_{e})^{n}$}\label{eq:cong of monge in thm ln k}\end{equation}
weakly as measures.

\begin{proof}
In the following proof it will be convenient to let $C$ denote a
sufficiently large constant (which may hence vary from line to line).
First observe that taking the logarithm of the inequality \ref{eq:exp decay}
in lemma \ref{lem:exponentia decay} immediately gives the upper bound
\[
k^{-1}\textrm{ln\,$K^{k}(z,z)\leq\phi_{e}(z)+Cn\ln k/k$}\]
To get a lower bound it is clearly enough to prove that for any fixed
point $z_{0}$ there is an element $f_{k}$ in $\mathcal{H}_{k}$
such that \begin{equation}
\left|f_{k}(z_{0})\right|^{2}e^{-k\phi_{e}(z_{0})}\geq1/C,\,\,\,\int_{\C^{n}}\left|f_{k}(z)\right|^{2}e^{-k\phi}\omega_{n}\leq C\label{pf of thm ln K: claim}\end{equation}
 where the constant $C$ is independent of $z_{0}$ and $k.$ To this
end we apply the Ohsawa-Takegoshi-Manivel extension theorem \cite{m}
(on the Stein manifold $\C^{n})$ with the continuous and \emph{strictly}
plurisubharmonic weight function $\psi_{k},$ where \begin{equation}
\psi_{k}(z):=k\phi_{e}(z)+C\ln(1+\left|z\right|^{2})\label{pf of thm ln k:weight fun}\end{equation}
This gives an extension of a constant, say $e^{k\phi_{e}(z_{0})/2},$
from the point $z_{0}$ to a holomorphic function $f_{k}$ in $\C^{n},$
such that \[
\left|f_{k}(z_{0})\right|^{2}e^{-k\phi_{e}(z_{0})}=1,\,\,\,\int_{\C^{n}}\left|f_{k}(z)\right|^{2}e^{-(k\phi_{e}+C\ln(1+\left|z\right|^{2}))}\omega_{n}\leq C\cdot1\]
 Note that, by a version of the Liouville theorem, this implies that
$f_{k}$ is a polynomial of degree at most $k.$ Finally, splitting
the latter integral with respect to the decomposition \[
\C^{n}=\left\{ \left|z\right|<C\right\} \bigcup\left\{ \left|z\right|\geq C\right\} ,\]
 then gives \ref{pf of thm ln K: claim}. Indeed, on $\left\{ \left|z\right|<C\right\} $
we have $\phi_{e}\leq\phi$ (by the definition \ref{eq:extem metric}
of $\phi_{e})$ and the factor $e^{-C\ln(1+\left|z\right|^{2})}$
is bounded by a constant. Furthermore, on $\left\{ \left|z\right|\geq C\right\} $
the assumption \ref{eq:ass on growth of phi} on the growth of $\phi$
implies that \[
k\phi_{e}+C\ln(1+\left|z\right|^{2})\leq k\phi\]
 All in all this shows that $\left\Vert f_{k}\right\Vert _{k\phi}^{2}\leq C,$
giving \ref{pf of thm ln K: claim}. 

The Monge-Ampere convergence \ref{eq:cong of monge in thm ln k} now
follows from the uniform convergence \ref{eq:conv of ln K in theorem ln K}
(see \cite{ko}).
\end{proof}
In fact, the uniform convergence in the previous theorem still holds
if $\phi$ is only assumed to be lower semi-continuous and locally
bounded. Indeed, then the use of the inequality $(i)$ in lemma \ref{lem:(Local-Morse-inequalities)}
may be replaced by the following Bernstein-Markov type inequality:
$B_{k}\leq C_{\epsilon}e^{k\epsilon}$ for any $\epsilon>0,$ whose
proof is an even simpler application of the submean property of holomorphic
functions. For a similar argument in the locally unbounded case see
\cite{berm3}. However, as shown in remark \ref{rem:holder} the $\mathcal{C}^{1,1}-$regularity
of $\phi$ is essential to obtain the optimal rate of convergence
$n\ln k/k$ in the theorem.

\subsection{The full Bergman kernel $K_{k}(z,w)$}

Combining the convergence in theorem \ref{thm:B in L1} with the local
inequalities \ref{eq:local morse on ball}, gives the following convergence
for the full Bergman kernel $K_{k}.$ The proof is completely analogous
to the proof of theorem 2.4 in part 1 of \cite{berm2}. 

\begin{thm}
\label{thm:k as meas}Let $K_{k}$ be the Bergman kernel of the Hilbert
space $\mathcal{H}_{k}$ of weighted polynomials in $\C^{n}$ of degree
less than $k.$ Then \[
\begin{array}{lr}
k^{-n}\left|K_{k}(z,w)\right|^{2}e^{-k\phi(z)}e^{-k\phi(w)}\omega_{n}(z)\wedge\omega_{n}(w)\rightarrow\Delta\wedge1_{D\cap X(0)}(dd^{c}\phi)^{n}/n!\end{array},\]
as measures on $\C^{n}\times\C^{n}$, in the weak {*}-topology, where
$\Delta$ is the current of integration along the diagonal in $\C^{n}\times\C^{n}.$
\end{thm}
Finally, we will show that around any interior point of the set $D\cap X(0)$
where $\phi$ is smooth the Bergman kernel $K_{k}(z,w)$ admits a
complete local asymptotic expansion in powers of $k,$ such that the
coefficients of the corresponding symbol expansion coincide with the
Tian-Zelditch-Catlin expansion for a positive Hermitian holomorphic
line bundle (see \cite{b-b-s} and the references therein). We will
use the notation $\phi(z,w)$ for a fixed almost holomorphic-anti-holomorphic
extension of $\phi$ from the diagonal $\Delta$ in $\C^{n}\times\C^{n},$
i.e. an extension such that the anti-holomorphic derivatives in $z$
and the holomorphic derivatives in $w$ vanish to infinite order on
$\Delta.$ 

\begin{thm}
\label{thm:asymp expansion} Let $K_{k}$ be the Bergman kernel of
the Hilbert space $\mathcal{H}_{k}$ of weighted polynomials in $\C^{n}$
of degree less than $k.$ Any interior point in $D\cap X(0)$ where
$\phi$ is smooth (i.e. in the class $\mathcal{C}^{\infty})$ has
a neighbourhood where $K_{k}(z,w)e^{-k\phi(z)/2}e^{-k\phi(w)/2}$
admits an asymptotic expansion as\begin{equation}
k^{n}(\det(dd^{c}\phi)(z)+b_{1}(z,w)k^{-1}+b_{2}(z,w)k^{-2}+...)e^{k\phi(z,w)},\label{eq:exp in prop}\end{equation}
where $b_{i}$ is a polynomial in the derivatives of $\phi,$ which
can be obtained by the recursion given in \cite{b-b-s}.
\end{thm}
\begin{proof}
Since the proof is essentially the same as the proof of theorem 4.4
in \cite{berm4} we will be very brief. One proceeds by adapting the
construction in \cite{b-b-s}, concerning \emph{postive} Hermitian
line bundles, to the present situation. The approach in \cite{b-b-s}
is to first construct a {}``local asymptotic Bergman kernel'' close
to any point where $\phi$ is smooth and $dd^{c}\phi>0$. Hence, the
local construction applies to the present situation as well. Then
the local kernel is shown to differ from the true kernel by a term
of order $O(k^{-\infty}),$ by solving a $\overline{\partial}$-equation
with a good $L^{2}-$estimate. In the present situation one applies
the $L^{2}-$estimates of Hörmander with the weight function $\psi_{k}$
(formula \ref{pf of thm ln k:weight fun}) occurring in the proof
of theorem \ref{thm:ln K}.
\end{proof}
Note that, even if the weight $\phi$ is smooth and globally psh such
an asymptotic expansion in integer powers of $k$ may not hold at
points where $\phi$ is not strictly psh (i.e. on $X-X(0)).$ Indeed,
at such a point $z_{0},$ say $z_{0}=0,$ the leading term (i.e. the
one of order $k^{n}$$)$ vanishes by lemma \ref{lem:(Local-Morse-inequalities)},
but taking $\phi(z)=\left|z\right|^{2+\epsilon},$ where $0<\epsilon<<1,$
gives $B_{k}(0)\geq Ck^{n/(1+\epsilon/2)},$ according to \ref{eq:holder ex}
applied to $\delta=-\epsilon.$

\section{\label{sec:Weighted-polynomials-associated}Weighted polynomials
associated to a Newton polytope}

In this section we briefly point out that the previous results extend
to the setting of weighted polynomials associated to a given Newton
polytope (which we show may even be taken to be a convex body). This
generalizes results of Shiffman-Zelditch \cite{sz4}. The main new
feature here compared to \cite{sz4} is that it is not assumed that
the weight function $\phi$ is plurisubharmonic nor invariant under
the action of the real unit torus (compare example \ref{exa:s-z}
for the setting in \cite{sz4}). 

Let $\Delta$ be a given convex body (i.e. compact with non-empty
interiour) in $\R^{n}$ and let $H_{\Delta}(z)$ be its support function,
considered as a function on the complex torus $\C^{*n}:$ \[
H_{\Delta}(z):=2\sup_{p\in\Delta}\ln(\left|z_{1}^{p_{1}}\right|\cdots\left|z_{n}^{p_{n}}\right|).\]
The role of the Lelong class is now played by the following class
of plurisubharmonic functions on $\C^{*n}:$ \[
\mathcal{L}_{\Delta}(\C^{*n}):=\left\{ \psi:\psi\,\textrm{psh,\,\,}\psi(z)\leq H_{\Delta}(z)+O(1)\right\} \]
and the role of the equilibrium potential $\phi_{e}$ by $\phi_{\Delta,e}$
defined by taking the sup over all $\widetilde{\phi}\leq\phi$ such
that $\widetilde{\phi}$ is in $\mathcal{L}_{\Delta}(\C^{*n}).$ Accordingly,
write $D_{\Delta}:=\{\phi_{\Delta,e}=\phi\}.$

Similarly, denote by $\mathcal{H}_{k\Delta}$ the Hilbert space spanned
by all Laurent monomials $z^{\alpha},$ where $\alpha$ is a vector
in $k\Delta\cap\Z$ (using multi index notation): \[
\mathcal{H}_{k\Delta}:=\bigoplus_{\alpha\in k\Delta\cap\Z^{n}}\C z^{\alpha}\]
with the norm \begin{equation}
\left\Vert f_{k}\right\Vert _{k\phi}^{2}:=\int_{\C^{*n}}\left|f_{k}(z)\right|^{2}e^{-k\phi}\omega_{n},\label{eq:scalar product}\end{equation}
assuming that $\phi$ is sufficiently large at the {}``boundary''
of $\C^{*n}$ so as to make the norms finite for $k$ large. More
precisely, after a {}``change of frame'' we may assume that zero
is an interiour point of $\Delta$ (just take $\textrm{exp }g(z)$
as in the beginning of the proof of lemma \ref{lem:(Local-Morse-inequalities)}
of the form $z^{\alpha_{0}}).$ Then the growth assumption is that
that there is an $\epsilon>0$ such that \begin{equation}
\phi(z)\geq(1+\epsilon)H_{\Delta}(z)\label{eq:growth ass newton}\end{equation}
outside some compact subset of $\C^{*n}.$ Finally, note that if $f_{k}$
is in $\mathcal{H}_{k\Delta},$ then $k^{-1}\ln\left|f_{k}(z)\right|^{2}$
is in $\mathcal{L}_{\Delta}(\C^{*n}).$ 

All previous results extend fairly straightforwardly to this setting.
For example, we have the following analog of theorem \ref{thm:B in L1}:

\begin{thm}
\label{thm:conv of B newton}Let $B_{k}$ be the Bergman function
of the Hilbert space $\mathcal{H}_{k\Delta}$ of weighted polynomials
in $\C^{n}$ associated to the Newton polytope $k\Delta.$ Then \begin{equation}
k^{-n}B_{k}\rightarrow1_{D_{\Delta}\cap X(0)}\det(dd^{c}\phi)\label{eq:l1 conv of B}\end{equation}
in $L^{1}(\C^{*n}).$ In particular, $k^{-n}B_{k}\omega_{n}$ conerges
weakly to the corresponding equilibrium measure $(dd^{c}\phi_{\Delta,e})^{n}/n!.$
\end{thm}
To see how to modify the previous proof first observe that the following
integrated versions of the theorem (that are used in its proof as
before), generalizing \ref{eq:mass bound in lelong} and \ref{eq:dim formel for B}
respectively, hold: \[
(i)\,\int_{\C^{*n}}(dd^{c}\phi_{e})^{n}/n!\leq\textrm{Vol $(\dot{\Delta})$},\,\,\,(ii)\,\dim\mathcal{H}_{k\Delta}\geq\textrm{Vol $(\dot{\Delta})$}k^{n}+o(k^{n}),\]
 where $\dot{\Delta}$ denotes the interiour of $\Delta.$ %
\footnote{In fact \emph{equalities} do hold. But the inequalities given are
the easiest ones to prove and then the equalities actually follow
from integrating the asymptotics in theorem \ref{thm:conv of B newton},
whose proof only uses the inequalies above.%
} To see that $(i)$ holds one first observes that, by a simple variant
of the monotonicity property \ref{eq:comp of mass}, it is enough
to prove $(i)$ with $\phi_{e}$ replaced by $H_{\Delta}(z).$ But
since $h_{\Delta}(v):=H_{\Delta}(e^{v_{1}},...e^{v_{n}})$ is convex
w.r.t $v\in\R^{n}$ $(i)$ then follows from well-known convex analysis,
using that the image of the subgradient $dh_{\Delta}$ is contained
in $\dot{\Delta}$ (cf. \cite{gr}). Finally, $(ii)$ is obtained
by approximating the number of lattice points in $\dot{k\Delta}$
by a Riemann sum in a standard fashion.

Finally, to prove theorem \ref{thm:conv of B newton} one also uses
that the analog of the regularity proposition \ref{pro:ref of equil pot}
holds for $\phi_{\Delta,e}$ on $\C^{*n},$ by replacing the action
of the translations $z\mapsto z+h$ by the action of the action of
the complex torus $\C^{*n},$ which may be represented by $z\mapsto e^{h}z$
(where $h\in\C^{n}).$ The rest of the proof of theorem \ref{thm:conv of B newton}
proceeds precisely as before.

\begin{example}
\label{exa:polyt}Let $\Delta\subset\Delta'$ where $\Delta'$ is
another convex body such that $\Delta$ is contained in its interiour.
Let $\phi$ be any smooth function such that $\phi=H_{\Delta'}+O(1).$
Then $\phi$ clearly satisfies the growth assumption \ref{eq:growth ass newton}.
If one further assumes that $\phi$ is invariant under the action
of the real unit torus (i.e. $\phi(z)=\Phi(v),$ where $v_{i}=\ln\left|z_{i}\right|^{2}$
and that $\phi$ is plurisubharmonic, so that $\Phi(v)$ is a convex
function on $\R^{n},$ then it is not hard to see that $D_{\Delta}=(d\Phi)^{-1}(\Delta).$
\end{example}
Specializing further gives the setting of Shiffman-Zelditch \cite{sz4}:

\begin{example}
\label{exa:s-z} Take $\Delta'$ to be the unit-simplex, $\Delta$
a rational polytope and $\phi(z)=\phi_{FS}(z):=\ln(1+\left|z\right|^{2})$
i.e. the potential of the Fubini-study metric. Note however that in
\cite{sz4} it is not assumed that $\Delta$ is contained in the \emph{interiour}
of the unit-simplex, which is related to the fact that $\phi_{FS}$
has the further property of extending smoothly to a metric on the
line bundle $\mathcal{O}(1)\rightarrow\P^{n}$ - compare the remark
below. Note also that specializing further to the case $n=1$ (i.e.
$\Delta'=[0,1]$ and $\Delta=[a,b]$ where $a,b\in\Q\bigcap]0,1[$
shows that the equlibrium weight $\phi_{\Delta,e}$ will in general
not be $\mathcal{C}^{1,1}$ on the complement of $\C^{*n}$ in $\C^{n},$
since $\phi_{\Delta,e}=H_{\Delta}(z)+O(1)=\ln(\left|z\right|^{2a}+\left|z\right|^{2b})+O(1),$
which tends to $-\infty$ as $z$ tends to zero.
\end{example}
\begin{rem}
\label{rem:boundary of pol}. Consider again the setting in the beginning
of example \ref{exa:polyt}, but assume also that $\Delta$ and $\Delta'$
both are rational polytopes. If $\Delta$ is to be allowed to intersect
the boundary of $\Delta'$ then one must assume that $\phi$ extends
to a metric on line bundle $\  L_{\Delta'}\rightarrow X_{\Delta'}$
over the corresponding toric variety $X_{\Delta}$ (compare section
5 in \cite{berm4} and references therein for notation). In this case
the results on Bergman kernel asymptotics etc follow from the results
announced in section 1.3 in \cite{berm4}, concerning \emph{subspace
versions} of the results in \cite{berm4}. The point is that $k\Delta$
defines a multiplier ideal sheaf $\mathcal{I}_{kH_{\Delta}}$ and
the global sections of the corresponding multiplier ideal sheaf $L_{\Delta'}^{\otimes k}\otimes\mathcal{I}_{kH_{\Delta}}\rightarrow X_{\Delta'}$
defines a sub Hilbert space $\mathcal{H}_{k\Delta}$ of the space
of all global sections $H^{0}(X_{\Delta'},L_{\Delta}^{\otimes k})$
equiped with the norm induced by the given metric $\phi$ on $L_{\Delta'}.$
\end{rem}

\end{document}